\documentclass[11pt]{amsart}

\usepackage{a4wide}
\usepackage{amsfonts}
\usepackage{psfrag}

\usepackage{amsmath}
\usepackage{amssymb}
\usepackage{amsthm}
\usepackage{graphicx}
\usepackage{color}

\newtheorem {definition}{Definition}

\newtheorem {thm}[section]{Theorem}
\newtheorem {rem}[definition]{Remark}

\newcommand{\cB}{{\mathcal B}}

\newcommand{\cD}{{\mathcal D}}

\newcommand{\cG}{{\mathcal G}}

\newcommand{\D}{\mathbb{D}}
\newcommand{\E}{\mathbb E}
\newcommand{\N}{\mathbb{N}}
\renewcommand{\P}{\mathbb{P}}
\newcommand{\EE}{\mathbb{E}}
\newcommand{\R}{\mathbb{R}}

\newcommand{\bg}{\mathbf{g}}
\newcommand{\bh}{\mathbf{h}}

\newcommand{\dd}{\mathsf{d}}
\newcommand{\ee}{\mathsf{e}}

\begin{document}
\title{On entire moments of self-similar Markov processes}

\author{M\'aty\'as Barczy}
\address{M\'aty\'as Barczy, University of Debrecen, Faculty of Informatics, Pf.12, H-4010 Debrecen, Hungary}
\email{barczy.matyas@inf.unideb.hu}

\author{Leif D\"oring}
\address{Leif D\"oring, Laboratoire de Probabilit\'es et Mod\'eles Al\'eatoires Universit\'e Paris 6, 4 place Jussieu, 75252 Paris }
\email{leif.doering@upmc.fr}
\thanks{
 M. Barczy has been {partially} supported by the Hungarian Scientific Research Fund under Grant No.\ OTKA T-079128, by the Hungarian
 Chinese Intergovernmental S \& T Cooperation Programme for 2011-2013 under Grant No.\ 10-1-2011-0079 and by the T\'AMOP-4.2.1/B-09/1/KONV-2010-0005 project.
 L. D\"oring has been supported by the Foundation Science Mat\'ematiques de Paris.
}



\begin{abstract}
It has been shown by Bertoin and Yor \cite{BY2} that {the law of} positive self-similar Markov processes (pssMps)
 that only jump downwards and do not hit zero in finite time
 are uniquely determined by their entire moments for which explicit formulas have been derived.
We use a recent jump-type stochastic differential equation approach to reprove and to extend their formulas.
\end{abstract}

\maketitle

\section*{Introduction and Results}
This article is focused on positive self-similar Markov processes (pssMps for short) introduced by Lamperti \cite{L} under
 the original name of semi-stable processes.
Let $\R_+$ denote the set of non-negative real numbers
 and let $\D$ be the space of c\`{a}dl\`{a}g functions $\omega:\R_+\to\R_+$ (right continuous with left limits)
 endowed with the Borel sigma-field $\cD$ generated by Skorokhod's $J_1$ topology.
Let us consider a family of probability measures $\{\P_z,z\geq 0\}$ on $(\D,\cD)$ under which the coordinate
 process $Z_t(\omega):=\omega(t)$, $t\geq 0$, is a strong Markov process starting from $z$ and fulfills
 the following scaling property: There exists a constant $a>0$, called the index of self-similarity, such that
 \begin{align}\label{self_sim}
   \text{the law of $(c^{-a}Z_{ct})_{t\geq 0}$ under $\P_z$ is $\P_{c^{-a}z}$}
 \end{align}
 for all $c>0$ and $z\geq 0$.
These laws are called non-negative self-similar Markov distributions and simultaneously the canonical process
 $Z$ on the probability space $(\D,\cD,\P_z)$ is called a non-negative self-similar Markov process (nnssMp for short)
 started from $z$.
Let us consider the process absorbed at the (possibly infinite) first hitting time of zero, $T_0$, given by
 \begin{align*}
	Z_t^\dag:=Z_t\mathbf 1_{\{t \leq T_0\}},\quad t\geq 0, \quad\text{with}\quad T_0:=\inf\{t\geq 0: Z_t=0\}.
 \end{align*}
The laws of the process $(Z_t^\dag)_{t\geq 0}$ under $\{\P_z,z\geq 0\}$, denoted by $\{\P_z^\dag, z\geq 0\}$,
 are non-negative self-similar Markov distributions with the same index of self-similarity.
They will be called positive self-similar Markov distributions and simultaneously, $(Z_t^\dag)_{t\geq 0}$ is called
 a positive self-similar Markov process.\\
By our convention, the difference between positive and non-negative self-similar Markov processes is only that pssMps have $0$ as a trap.
Since the index of self-similarity can always be transformed to $1$ by taking the power $1/a$ of $Z$, in what follows we can assume without loss of generality that $a=1$.

Lamperti \cite{L} constructed all pssMps as time-changed exponentials of L\'evy processes; more precisely, if $\xi$ is a
 (possibly killed, extended real-valued) L\'evy process, he showed that
 \begin{align*}
    Z_t:= z \exp\left(\xi_{\tau(tz^{-1})}\right), \qquad 0\leq t<T_0,
 \end{align*}
 where
 \[
   \tau(t):=\inf\{s\geq 0 : I_s\geq t\} \qquad \text{and} \qquad I_t:=\int_0^t \exp(\xi_s)\,\dd s,\qquad t\geq 0,
 \]
 is a pssMp of self-similarity index $1$ started from $z$.
Conversely, any pssMp of self-similarity index $1$ can be represented in this way for some L\'evy process $\xi$ possibly killed at an independent exponential time $\zeta$.
In case of killing the coffin state is chosen to be $-\infty$.
A remarkable consequence of Lamperti's representation is that a pssMp
 started at $z>0$ hits zero in finite time almost surely if and only if
 $\xi$ has infinite lifetime (i.e., no killing) and drifts to $-\infty$
 (i.e., $\lim_{t\to\infty}\xi_t=-\infty$ almost surely) or $\xi$ has finite lifetime,
 see Lamperti \cite[Lemmas 2.5 and 3.2]{L}.
In what follows we consider pssMps that only jump downwards, which is equivalent to imposing the assumption
\begin{align*}
	\textbf{(A1) }\quad \xi \text{ is spectrally negative}.
\end{align*}
In recent years, several authors extended Lamperti's characterization to nnssMps. Under Assumption \textbf{(A1)}, Bertoin and Yor \cite{BY2}
showed that pssMps that do not hit zero in finite time  (i.e. the Lamperti transformed L\'evy process $\xi$ is neither killed nor  drifts to $-\infty$)
can be extended uniquely to a Feller process on $[0,\infty)$, and {the law of} the corresponding nnssMps that start from zero non-trivially were fully characterized by their entire moments.
Much deeper results for general pssMps that do not hit zero, and not necessarily satisfy Assumption \textbf{(A1)}, have been obtained by
Caballero and Chaumont \cite{CC} and Bertoin and Savov \cite{BS}. For pssMps that do hit zero in finite time almost surely, Rivero \cite{R1}, \cite{R2} and
Fitzsimmons \cite{F} proved independently the existence and uniqueness of a recurrent extension that leaves zero continuously if and only if the Cram\'er type condition
  \begin{align}\label{cramer}
	{\text{there is a }0<\theta<1\text{ such that }\Psi(\theta)=0 }
 \end{align}
 holds, where $\Psi(\lambda):=\log \E(\ee^{\lambda \xi_1};\zeta>1)$ denotes the Laplace exponent of $\xi$
 killed at $\zeta$  (we simply used $\E$ instead of $\E_z$, since the law of $\xi$ does not depend on $z$).
Note that under Assumption \textbf{(A1)} the Laplace exponent $\lambda\mapsto \Psi(\lambda)$ is a convex, continuous function
 on $\R_+$ with derivative $\E(\xi_1)$ at zero, if $\xi_1$ has finite mean, and hence, if additionally $\xi$ drifts to $-\infty$
 or $\xi$ being killed, then it follows readily that (see, e.g. Kyprianou \cite[page 81]{K}) the assumption
 \begin{align*}
	\textbf{(A2) }\quad \Psi(1)>0
 \end{align*}
 holds if and only if Condition \eqref{cramer} is satisfied, i.e. precisely in the settings of Fitzsimmons \cite{F}
 and Rivero \cite{R1}, \cite{R2}.
We also note that a companion result of Bertoin and Yor \cite[Proposition 2]{BY2} has been obtained by Patie \cite[Theorem 2.3]{Pat}
 who characterized the law of the exponential functional $I_\infty:=\lim_{t\to\infty}I_t$ by an explicit form of its Laplace transform
 provided that Assumptions \textbf{(A1)} and \textbf{(A2)} hold and $\xi$ drifts to $-\infty$.
\smallskip

The aim of this article is to use a jump-type SDE characterization of nnssMps given in D\"oring and Barczy \cite{DB} to prove the entire moment formulas of Bertoin and Yor \cite{BY2} for nnssMps $Z$ for which $Z^\dag$ satisfies \textbf{(A1)} and \textbf{(A2)}.
We call the attention that, in contrast to the other known extensions of Lamperti's characterization to zero initial condition $z=0$,
 the jump-type SDE approach handles the cases $\xi$ drifting to $-\infty$, $\xi$ oscillating, $\xi$ drifting to $+\infty$ or
 $\xi$ being killed at once.

\begin{thm}\label{Prop:Bertoin_Yor}
Let $\{\P_z, z\geq 0\}$ be a {family of non-negative self-similar Markov distributions} of self-similarity index 1 such that
 the Lamperti transformed L\'evy process for the corresponding family $\{\P_z^\dag, z\geq 0\}$ of positive self-similar Markov distributions
 satisfies the Assumptions \textbf{(A1)} and \textbf{(A2)} with Laplace exponent $\Psi$. Then the {laws}  $\{\P_z, z\geq 0\}$
 {are} uniquely determined by the entire moment formulas
 \begin{align}\label{Moments}
    \EE_z(Z_t^n) = z^n + \sum_{\ell=1}^n \frac{\Psi(n)\cdots \Psi(n-\ell+1)}{\ell!} z^{n-\ell}t^\ell,  \qquad n\in\N, t\geq 0, z\geq 0.
 \end{align}
\end{thm}
Let us quickly recall how Bertoin and Yor \cite{BY2} proved the moment formulas when $\xi$ does not drift to $-\infty$.
Since zero is never hit, the full information on $Z$ started at $z>0$ is already contained in the infinitesimal generator which was determined by Lamperti \cite[Theorem 6.1]{L}.
Applying the infinitesimal generator to the function $z\mapsto z^n$  with  $n>0$, Proposition VII.1.2 in Revuz and Yor \cite{RY} yields the recursion
 \[
   \frac{\partial}{\partial t} \EE_z\big(Z_t^n\big) = \Psi(n)\, \EE_z\big(Z_t^{n-1}\big),\qquad t\geq0,\,n\in\N.
 \]
Iterating this equation one can derive formula \eqref{Moments}.
They extended the formula to $z=0$ by sending $z$ to zero and by showing that the moment problem for the law of $Z_t$
  under $\P_z$, $z\geq0$, is well-posed.\\
The proof presented here is based on a reformulation of Lamperti's infinitesimal generator by jump type stochastic differential equations. The striking feature of the jump type SDEs is that they can be readily extended after hitting zero and, thus, entire moment formulas for recurrent extensions can be derived from It\=o's formula also if Lamperti's generator characterization is not available.
\begin{rem}
Bertoin and Yor \cite{BY2} {determined the law of} the nnssMps started from zero via the moment formulas \eqref{Moments}. Conversely, we use our construction from D\"oring and Barczy \cite{DB} of the {nnssMps} started from zero to derive the moment formulas.
It is not clear to us how the moment formulas can be deduced directly from Lamperti's transformation in the general case.
\end{rem}

\section*{Proof of Theorem \ref{Prop:Bertoin_Yor}}

From now on we will work on a stochastic basis $(\Omega,\cG,(\cG_t)_{t\geq 0},P)$ satisfying
 the usual conditions and we will denote the underlying probability measure
 and the expectation with respect to it by $P$ and $E$, respectively, instead of $\P_z$ and $\EE_z$
 (in contrast to the introduction, but without confusion).

Suppose that $(Z_t)_{t\geq 0}$ is a nnssMp and $(Z^\dag_t)_{t\geq 0}$ is the corresponding pssMp trapped when first hitting zero.
Since $Z^\dag$ is a pssMp there is a L\'evy process $\xi$ (possibly killed at rate $q$) with Laplace exponent $\Psi$
 and L\'evy triplet $({\gamma},\sigma^2,\Pi)$ corresponding to $Z^\dag$ under Lamperti's representation.
{Recall that $\gamma\in\R$, $\sigma\geq 0$ and $\Pi$ is a deterministic measure on $(-\infty,0)$ satisfying
 $\int_{-\infty}^0\min(1,u^2)\,\Pi(\dd u)<\infty$.}
In D\"oring and Barczy \cite{DB} it was shown that
 for all $z\geq 0$ the law of $Z$ coincides with the law of the pathwise unique strong solution to
 \begin{align}\label{SDE_pssMp}
   \begin{split}		
	Z_t&=z+\Psi(1)\,t + \sigma\int_0^{t } \sqrt{Z_s} \dd B_s\\
           &\quad+\int_0^{t }\int_0^{\infty}\int_{-\infty}^0 \mathbf 1_{\{r Z_{s-}\leq 1\}}Z_{s-}[\ee^u-1](\mathcal N-\mathcal N')(\dd s,\dd r,\dd u)\\
           &\quad-\int_0^t\int_0^\infty\mathbf 1_{\{r Z_{s-}\leq 1\}} Z_{s-}(\mathcal M-\mathcal M')(\dd s,\dd r),
           \qquad t\geq 0,
   \end{split}
 \end{align}
 where $B$ is a standard Wiener process, $\mathcal N$ is a Poisson random measure
 on $(0,\infty)\times (0,\infty)\times (-\infty,0)$ with intensity measure $\mathcal N'(\dd s,\dd r,\dd u)=\dd s\otimes \dd r \otimes \,\Pi(\dd u)$,
 $\mathcal M$ is a Poisson random measure on $(0,\infty)\times (0,\infty)$ with intensity measure $\mathcal M'(\dd s,\dd r)=q\,\dd s\otimes \dd r$
 and $q\geq 0$ is the killing rate (i.e. $P(\zeta>1)=\ee^{-q}$).
If we abbreviate
 \begin{align*}
  \bg(x,r,u):=\mathbf 1_{\{ rx\leq 1\}} x (\ee^u-1)\quad\text{ and }\quad
    \bh(x,r):=-\mathbf 1_{\{ rx\leq 1\}} x,
 \end{align*}
 than It\={o}'s formula for non-continuous semi-martingales (see, e.g. Di Nunno et al. \cite[Theorem 9.5]{NOP}
 or Ikeda and Watanabe \cite[Chapter II, Theorem 5.1]{IW}) implies, for $t\geq 0$,
 \begin{align*}
    Z_t^n   &= z^n + \int_0^t n Z_s^{n-1}\sigma \sqrt{Z_s}\,\dd B_s  +\Psi(1) \int_0^t n Z_s^{n-1}\,\dd s\\
            &\quad  + \frac{1}{2}\int_0^t n(n-1)Z_s^{n-2}\sigma^2 Z_s\,\dd s\\
            &\quad  + \int_0^{t }\int_0^{\infty}\int_{-\infty}^0\Big[ \big(Z_s + \bg(Z_s,r,u) \big)^n- Z_s^n
                     - \bg(Z_s,r,u) n Z_s^{n-1}\Big]\dd s\,\dd r\,\Pi(\dd u) \\
            &\quad  + \int_0^{t }\int_0^{\infty}\int_{-\infty}^0\Big[ \big(Z_{s-} + \bg(Z_{s-},r,u) \big)^n- Z_{s-}^n \Big] (\mathcal N-\mathcal N')(\dd s,\dd r,\dd u) \\
            &\quad  + q\int_0^{t }\int_0^{\infty}\Big[ \big(Z_s + \bh(Z_s,r) \big)^n- Z_s^n
                               - \bh(Z_s,r) n Z_s^{n-1}\Big]\dd s\,\dd r \\
            &\quad  + \int_0^{t }\int_0^{\infty}\Big[ \big(Z_{s-} + \bh(Z_{s-},r) \big)^n- Z_{s-}^n \Big] (\mathcal M-\mathcal M')(\dd s,\dd r)
 \end{align*}
 which, carrying out the compensator integrals, gives
  \begin{align*}
   Z_t^n &=z^n + n\sigma\int_0^t Z_s^{n-\frac{1}{2}}\,\dd B_s \\
         &\quad +\left[n\Psi(1)+\frac{n(n-1)\sigma^2}{2}+\,\int_{-\infty}^0 \big[ \ee^{nu} - 1 - n(\ee^u-1) \big]\Pi(\dd u)+q(n-1)\right]
                      \int_0^t Z_s^{n-1}\,\dd s\\
         &\quad + \int_0^t\int_0^{\infty}\int_{-\infty}^0 \mathbf 1_{\{ rZ_{s-} \leq 1\}} (\ee^{nu} - 1) Z_{s-}^n (\mathcal N-\mathcal N')(\dd s,\dd r,\dd u)  \\
         &\quad - \int_0^t\int_0^{\infty} \mathbf 1_{\{ rZ_{s-} \leq 1\}}  Z_{s-}^n (\mathcal M-\mathcal M')(\dd s,\dd r).
 \end{align*}
By Sato \cite[Theorem 25.17]{Sat}, we see that
\begin{align*}
	&\quad n\Psi(1)+\frac{n(n-1)\sigma^2}{2}+\int_{-\infty}^0 \big[(\ee^{nu} - 1 - n(\ee^u-1) \big]\Pi(\dd u)+q(n-1) \\
	& =n\gamma+\frac{n\sigma^2}{2}+n\int_{-\infty}^0(\ee^u-1-u\textbf{1}_{\{|u|\leq 1\}})\Pi(\dd u) - nq
        + \frac{n(n-1)\sigma^2}{2}\\
	&\quad+\int_{-\infty}^0 \big[ (\ee^{nu} - 1 - n(\ee^u-1) \big]\Pi(\dd u) +q(n-1)\\
	&=n\gamma+\frac{n^2\sigma^2}{2}	+ \int_{-\infty}^0 \big[ (\ee^{nu}-1-n  u\textbf{1}_{\{|u|\leq1\}}) \big]\Pi(\dd u) - q\\
	&= \log E[\ee^{n\xi_1};\zeta>1] =\Psi(n).
\end{align*}
Plugging-in, we obtain
 \begin{align}\label{SDE_pssMp_moment}
   \begin{split}
    Z_t^n &= z^n +\Psi(n) \int_0^t Z_s^{n-1}\,\dd s+ n\sigma\int_0^t Z_s^{n-\frac{1}{2}}\,\dd B_s\\
             &\quad + \int_0^t\int_0^{\infty}\int_{-\infty}^0\mathbf 1_{\{ rZ_{s-} \leq 1\}}(\ee^{nu} -1) Z_{s-}^n (\mathcal N-\mathcal N')(\dd s,\dd r,\dd u)\\
                          &\quad - \int_0^t\int_0^{\infty}\mathbf 1_{\{ rZ_{s-} \leq 1\}}Z_{s-}^n (\mathcal M-\mathcal M')(\dd s,\dd r)\\
		 &=: z^n +\Psi(n) \int_0^t Z_s^{n-1}\,\dd s+M^{(1)}_t+M^{(2)}_t+M^{(3)}_t,\qquad t\geq 0.
   \end{split}
 \end{align}
In what follows we want to take expectations to deduce the recursive equations
\begin{align}\label{Recursive}
	E\big(Z_t^n\big)=z^n +\Psi(n) \int_0^t E\big(Z_s^{n-1}\big)\,\dd s,\quad n\in\N,\; t\geq 0.
\end{align}
Before doing so, we need to show that $(M^{(1)}_t)_{t\geq 0}, (M^{(2)}_t)_{t\geq 0}$ and $(M^{(3)}_t)_{t\geq 0}$ are
 martingales and not only local martingales with respect to the filtration $(\cG_t)_{t\geq 0}$.
First we check that they are local martingales which easily follows by the construction of the stochastic integrals.
Indeed, let $\delta_m:=\inf\{ t\geq 0 : Z_t\geq m\}$, $m\in\N$.
 Then, by Step 1a in the proof of Proposition 3.13 in D\"oring and Barczy \cite{DB},
 we have \ $P(\lim_{m\to\infty} \delta_m = \infty)=1$ and using $(\delta_m)_{m\in\N}$ as a localizing sequence,
  for $M^{(1)}$ we have
  \begin{align*}
   E\left( \int_0^{t\wedge\delta_m} Z_s^{2n-1}\,\dd s \right)
      \leq m^{2n-1}E(t\wedge \delta_m) \leq t m^{2n-1},
  \end{align*}
 for $M^{(2)}$,
  \begin{align*}
 &\quad  E\left( \int_0^{t\wedge \delta_m}\int_0^{\infty}\int_{-\infty}^0
         \mathbf 1_{\{ rZ_{s} \leq 1\}}
           (\ee^{nu} -1)^2 Z_{s}^{2n}\dd s\,\dd r\,\Pi(\dd u)
       \right)\\
     & = E\left(\int_0^{t\wedge \delta_m}Z_{s}^{2n-1}\,\dd s\right)
        \int_{-\infty}^0 (\ee^{nu} - 1)^2\,\Pi(\dd u)\\
     & \leq m^{2n-1} E(t\wedge \delta_m) \int_{-\infty}^0 (\ee^{nu} - 1)^2 \Pi(\dd u)
      <\infty,
 \end{align*}
 and for $M^{(3)}$,
 \begin{align*}
&\quad  E\left( \int_0^{t\wedge \delta_m}\!\!\!\int_0^\infty
         \mathbf 1_{\{ rZ_{s} \leq 1\}} Z_{s}^{2n} q \dd s\,\dd r \right)\\
   & = q E\left( \int_0^{t\wedge \delta_m}Z_{s}^{2n-1} \dd s \right)
    \leq q m^{2n-1} E(t\wedge \delta_m) <\infty.
 \end{align*}
In case $M^{(2)}$ for the last step we used the asymptotic equivalence $(\ee^{nu}-1)^2\sim n^2u^2$ at zero and
 the integrability property $\int_{-\infty}^0 \min(1,u^2)\Pi(\dd u)<\infty$ for the L\'evy measure $\Pi$.
The desired local martingale property of $M^{(i)}$, $i=1,2,3$, follows by Ikeda and Watanabe \cite[pages 57 and 62]{IW}.
Hence, taking expectations in the SDE \eqref{SDE_pssMp_moment}, we find that
 \begin{align*}
   E(Z_{t\wedge \delta_m}^n) &= z^n + \Psi(n)E\bigg(\int_0^{t\wedge \delta_m}Z_s^{n-1}\,\dd s\bigg)
     \leq  z^n + \Psi(n)\int_0^{t} E\big(Z_{s\wedge \delta_m}^{n-1}\big)\,\dd s
 \end{align*}
 for all $t\geq 0$ and $m\in\N$.
Here the inequality follows by the following two facts:
 \begin{itemize}
   \item
         \[
            \int_0^{t\wedge \delta_m}Z_s^{n-1}\,\dd s
             \leq \int_0^t Z_{s\wedge \delta_m}^{n-1}\,\dd s,\qquad t\geq 0, \; m\in\N,
         \]
         where the inequality (not being necessarily an equality) is explained by that for $t\geq\delta_m$, the left-hand side is $\int_0^{\delta_m} Z_s^{n-1}\,\dd s$
         and the right-hand side is $\int_0^{\delta_m} Z_s^{n-1}\,\dd s + (t-\delta_m)Z_{\delta_m}^{n-1}$.
    \item the assumption $\Psi(1)>0$ yields that $\Psi(n)>0$, $n\in\N$, due to the convexity of $\Psi$.
 \end{itemize}
By induction we obtain
 \begin{align*}
     E(Z_{t\wedge \delta_m}^n) \leq z^n +\sum_{\ell=1}^n \frac{\Psi(n)\cdots\Psi(n-\ell+1)}{\ell!}z^{n-\ell}t^\ell,
     \qquad t\geq 0,\,m\in\N,\,n\in\N.
 \end{align*}
Hence, by Fatou's lemma,
 \begin{align}\label{help_Fatou}
     E(Z_t^n) \leq \liminf_{m\to\infty}E(Z_{t\wedge \delta_m}^n) \leq z^n +\sum_{\ell=1}^n \frac{\Psi(n)\cdots\Psi(n-\ell+1)}{\ell!}z^{n-\ell}t^\ell,\qquad t\geq 0,\, n\in\N.
 \end{align}
Next, we can deduce that $M^{(1)}, M^{(2)}$ and $M^{(3)}$ are indeed true martingales.
First, we note that a local martingale $M$ is a martingale if $E(\sup_{t\in[0,T]}\vert M_t\vert)<\infty$ for all $T>0$
 (see, e.g. Theorem I.51 in Protter \cite{P}).
Since the finiteness of the second moment implies the finiteness of the first moment, the condition that $E(\sup_{t\in[0,T]} M_t^2)<\infty$ for all $T>0$
 also yields the martingale property of $M$.
First we give the argument for $M^{(1)}$.
By the Burkholder-Davis-Gundy inequality (see, e.g. Karatzas and Shreve \cite[Chapter 3, Theorem 3.28]{KS}) there exists a universal constant $C>0$ such that
 \[
   E\big( \sup_{t\in[0,T]} (M^{(1)}_t)^2\big) \leq C E\big(\langle M^{(1)}\rangle_T\big),
    \qquad \forall \, T>0,
 \]
 where $\langle\cdot\rangle$ denotes the quadratic variation process of a continuous local martingale.
Then
 \begin{align*}
   E\big( \sup_{t\in[0,T]} (M^{(1)}_t)^2\big) \leq Cn^2\sigma^2 E\left( \int_0^T Z_s^{2n-1}\,\dd s\right)<\infty,
 \end{align*}
 where the last inequality follows from Fubini's theorem and \eqref{help_Fatou}.
The same argument applies to $M^{(2)}$ and $M^{(3)}$ using the non-continuous version of the Burkholder-Davis-Gundy inequality (see, e.g. Dellacherie and Meyer \cite[Theorem VII. 92]{DM}) with the bracket process $[\cdot]$ for non-continuous semimartingales. To obtain the estimate, we use
 \begin{align*}
     E\big([ M^{(2)} ]_T\big)
     & = E\left( \int_0^T\int_0^{\infty}\int_{-\infty}^0 \mathbf 1_{\{ rZ_{s} \leq 1\}}(\ee^{nu} -1)^2 Z_{s}^{2n}
            \,\dd s\,\dd r\,\Pi(\dd u)\right)\\
     & = \int_0^T E(Z_{s}^{2n-1})\,\dd s \int_{-\infty}^0 (\ee^{nu} -1)^2  \Pi(\dd u) <\infty,
 \end{align*}
where the first equality follows by Jacod and Shiryaev \cite[Theorem I.1.33]{JS} and the last inequality by \eqref{help_Fatou} and the integrability property of $\Pi$.
 The finiteness of $E([ M^{(3)} ]_T)$ can be checked in the same way.
Thus, we verified that we can take the expectation of \eqref{SDE_pssMp_moment} to deduce \eqref{Recursive}.
Iterating Equation \eqref{Recursive} (which is in fact the same recursion for $E(Z_t^n)$ that was obtained by Bertoin and Yor \cite[page 38]{BY2})
 we find inductively
 \[
   E(Z_t^n) = z^n + \sum_{\ell=1}^n \frac{\Psi(n)\cdots \Psi(n-\ell+1)}{\ell!} z^{n-\ell} t^\ell,
                   \qquad \,n\in\N,t\geq 0,z\geq 0,
 \]
 as desired.
To deduce that, for all $t\geq 0$ and $z\geq 0$, the entire moments {\eqref{Moments}} of $Z_t$ uniquely determine its law,
 we only need to verify $E[\ee^{\theta Z_t}]<\infty$ for some $\theta>0$ sufficiently small. But this follows by expanding the exponential,
 plugging-in \eqref{Moments} and using that $\Psi(n)\leq K n^2$, $n\in\N$, with some positive constant $K>0$ for spectrally negative L\'evy processes.
For more details see the proof of part (ii) of Proposition 1 in Bertoin and Yor \cite{BY2}.

{Finally, using that $Z$ is a time-homogeneous (strong) Markov process (being the pathwise unique strong solution of the SDE \eqref{SDE_pssMp})
 and that the entire moments \eqref{Moments} determine uniquely the conditional law of the one-dimensional marginals of $Z$ given
 the initial value $Z_0$, we have the entire moments \eqref{Moments} determine uniquely the conditional finite dimensional marginals of $Z$
 (given the initial value $Z_0$) too.
Since $\{ \pi^{-1}_{t_1,\ldots,t_k}(B): \,\, 0\leq t_1\leq \cdots\leq t_k, \; B\in\cB(\R^k),\; k\in\N \}$
 is a separating class for $\D$, where $\pi_{t_1,\ldots,t_k}$ denotes the natural projection from $\D$ to $\R^k$,
 we have entire moments \eqref{Moments} determine uniquely the laws $\{\P_z,z\geq 0\}$ too, completing the proof.}

\section*{Acknowledgement}
We are grateful to the referee for several valuable
comments that have led to an improvement of the manuscript.


\begin{thebibliography}{99}

\bibitem{BS}
J. Bertoin and M. Savov:
 Some applications of duality for L\'evy processes in a half-line.
 {\it Bull. Lond. Math. Soc.} {\bf 43} 97--110, (2011).

\bibitem{BY2}
J. Bertoin and M. Yor:
On the entire moments of self-similar Markov processes and exponential functionals of L\'evy processes.
 {\it Ann. Fac. Sci. Toulouse Math.} {\bf S\'erie 6, Vol. 11} 33--45, (2002).

\bibitem{CC}
M.-E. Caballero and L. Chaumont:
Weak convergence of positive self-similar Markov processes and overshoots of L\'evy processes.
 {\it Ann. Probab.} {\bf 34} 1012--1034, (2006).

\bibitem{DM}
C. Dellacherie and P. A. Meyer:
 {\sl Probabilit\'es et potentiel. Chapitres V \`{a} VIII. Th\'eorie des martingales.}
Hermann, Paris, 1980.

\bibitem{NOP}
G. Di Nunno, B. {\O}ksendal and F. Proske: {\sl Malliavin calculus for L\'evy processes with applications to finance.
 Corrected second printing.} Springer, 2009.

\bibitem{DB}
L. D\"oring and M. Barczy:
A jump type SDE approach to positive self-similar Markov processes.
 {\it ArXiv} 1111.3235, (2011). URL: http://arxiv.org/abs/1111.3235

\bibitem{F}  P. Fitzsimmons: On the existence of recurrent extensions of self-similar
 Markov processes. {\it Electron. Comm. Probab.} {\bf 11} 230-241, (2006).

\bibitem{IW}
N. Ikeda and S. Watanabe: {\sl Stochastic differential equations and diffusion processes.}
 North-Holland Publishing Company, 1981.

\bibitem{JS}
J. Jacod and A.N. Shiryaev: {\sl Limit theorems for stochastic processes. Second edition.}
 Springer-Verlag, Berlin, 2003.

\bibitem{KS}
I. Karatzas and S.E. Shreve: {\sl Brownian motion and stochastic calculus. Second edition.}
 Springer-Verlag, 1991.

\bibitem{K}
A.E. Kyprianou: {\sl Introductory lectures on fluctuations of L\'evy process with applications.}
 Universitext, Springer, 2006.

\bibitem{L}
J. Lamperti:
Semi-stable Markov processes. I.
{\it Z. Wahr. und Verw. Gebiete} {\bf 22} 205--225, (1972).

\bibitem{Pat}
P. Patie: Infinite divisibility of solutions to some self-similar integro-differential equations and exponential
 functionals of L\'evy processes.
{\it Ann. Inst. H. Poincaré Probab. Statist.} {\bf 45}(3) 667--684, (2009).

\bibitem{P}
P. Protter:
{\sl Stochastic integration and differential equations. Second edition.}
Springer-Verlag, 2004.

\bibitem{RY} D. Revuz and M. Yor:
 {\sl Continuous martingales and Brownian motion. Third edition.}
 Springer-Verlag Berlin Heidelberg, 1999.

\bibitem{R1}
V. Rivero: Recurrent extensions of self-similar Markov processes and Cram\'er's condition.
 {\it Bernoulli} {\bf 11}(3) 471--509 (2005).

\bibitem{R2}
V. Rivero: Recurrent extensions of self-similar Markov processes and Cram\'er's condition. II.
 {\it Bernoulli} {\bf 13}(4) 1053--1070 (2007).

\bibitem{Sat}
K.-I. Sato: {\sl L\'evy processes and infinitely divisible distributions.}
 Cambridge University Press, Cambridge, 1999.

\end{thebibliography}
\end{document}